\documentclass[12pt]{article}

\usepackage[dvips]{graphicx}
\usepackage{t1enc,amssymb,amsbsy,amsthm,amsmath,amsfonts}

\renewcommand{\Re}{{\rm I\kern-0.16em R}}

\def\@begintheorem#1#2{\trivlist \item[\hskip \labelsep{\bf #1\ #2}]}
\def\@opargbegintheorem#1#2#3{\trivlist
      \item[\hskip \labelsep{\bf #1\ #2\ (#3)}]}

\newtheorem{proposition}{Proposition}

\newtheorem{thm}[proposition]{Theorem}

\def\P{{\bf P}}

\def\E{{\bf E}}

\def\cF{{\cal F}}

\def\al{\alpha}

\def\si{\sigma}

\begin{document}

\author{Paavo Salminen\\{\small Åbo Akademi University,}
\\{\small Mathematical Department,}
\\{\small Fänriksgatan 3 B,}
\\{\small FIN-20500 Åbo, Finland,} 
\\{\small email: phsalmin@abo.fi}
\and
Marc Yor\\
{\small Universit\'e Pierre et Marie Curie,}\\
{\small Laboratoire de Probabilit\'es }\\
{\small et Mod\`eles al\'eatoires ,}\\
{\small 4, Place Jussieu, Case 188}\\
{\small F-75252 Paris Cedex 05, France}
}
\vskip5cm

\title{
A note on a.s. finiteness of perpetual integral functionals of diffusions}
\date{}

\maketitle

\begin{abstract}
In this note, with the help of the boundary classification of diffusions, we derive 
a criterion of the convergence of  perpetual integral functionals of transient real-valued diffusions.

In the particular case of transient Bessel processes, we note that this criterion 
agrees with the one obtained via Jeulin's convergence lemma.  
\end{abstract}

{\bf 1.} Consider a diffusion $Y$ on an open interval $I=(l,r)$
determined by the SDE
$$
dY_t=\sigma(Y_t)\,dW_t+b(Y_t)\,dt,
$$
where $W$ is a standard Wiener process defined in a complete
probability space $(\Omega,\cF,\{\cF_t\},\P).$
It is assumed that $\sigma$ and $b$ are continuous and
$\si(x)>0$ for all $x\in I$. We assume also that $Y$ is transient  
and 
\begin{equation}
\label{a10}
\lim_{t\to\zeta} Y_t=r\quad {\rm a.s.},
\end{equation}
where $\zeta$ is the life time of
$Y.$ Hence, if $\zeta<\infty$ then 
$$
\zeta=H_r(Y):=\inf\{t:\, Y_t=r\}.
$$
For the speed and the scale measure of $Y$ we use 
\begin{equation}
\label{a11}
m^Y(dx)=2\,\sigma^2(x){\rm e}^{B^Y(x)}\,dx
\quad{\rm and}\quad 
S^Y(dx)={\rm e}^{-B^Y(x)}\,dx,
\end{equation}
respectively, where 
\begin{equation}
\label{B1}
B^Y(x)=2\,\int^x \frac{b(z)}{\sigma^2(z)}\,dz.
\end{equation}

Let $f$ be a positive and continuous function defined on $I,$ and
consider the perpetual integral functional 
$$
A_\zeta(f):=\int_0^\zeta f(Y_s)\,ds.
$$
We are interested in finding neccessary and sufficient conditions 
for a.s. finiteness of $A_\zeta(f).$
When $Y$ is a Brownian motion with drift $\mu>0$ such a condition is 
that the function $f$ is
integrable at $+\infty$ (see Engelbert and Senf \cite{engelbertsenf91} and 
Salminen and Yor \cite{salminenyor03}). This condition
is derived in \cite{salminenyor03} via Ray-Knight theorems and the
stationarity property of the
local time processes (which makes Jeulin's lemma applicable). 
In this note a condition (see Theorem \ref{!!}) valid for general $Y$ is deduced by exploiting the
fact that $A_{H_x}(f)$ for $x<r$ can, via random time change, be seen as the
first hitting time of a point for another diffusion.

{\bf 2.} The next proposition presents the key result connecting perpetual integral
functionals to first hitting times. The result is a generalization of 
Proposition 2.1 in \cite{salminenyor04} discussed in 
Propositions 2.1 and 2.3 in
\cite{borodinsalminen04}.

\begin{proposition}
\label{!}
 {\sl
Let $Y$ and $f$ be as above, and assume that there exists a two
times 
continuously differentiable function $g$ such that 
\begin{equation}
\label{fg}
f(x)=\big(g'(x)\si(x)\big)^2,\quad x\in I.
\end{equation}
Set for $t>0$ 
\begin{equation}
\label{A2}
A_t:=\int_0^t f(Y_s)\,ds.
\end{equation}
and let 
$\{a_t\,:\,0\leq t<A_\zeta\}$ denote the inverse of $A,$ that is,
$$
a_t:=\min\big\{s: A_s>t\big\}, \qquad t\in [0,A_\zeta).
$$
Then the process $Z$ given by
$$
Z_t:=g\left(Y_{a_t}\right), \qquad t\in[0,A_\zeta),
$$
is a diffusion satisfying the SDE
$$
dZ_t=d\widetilde W_t+G(g^{-1}(Z_t))\,dt, \qquad t\in[0,A_\zeta).
$$
where $\widetilde W_t$ is a Brownian motion and
$$
G(x)=\frac 1{f(x)}
\left( \frac 12\, \si(x)^2\, g''(x)
+b(x)\,g'(x)\right).
$$
Moreover, for $l<x<y<r$ 
\begin{equation}
\label{aa10}
A_{H_y(Y)}=\inf\{t:\, Z_t=g(y)\}=: H_{g(y)}(Z)\quad a.s.
\end{equation}
with $Y_0=x$ and $Z_0=g(x).$
}
\end{proposition}

{\bf 3.}  To fix ideas, assume that the function $g$ as introduced in
Proposition~\ref{!} is increasing. We define $g(r):=\lim_{x\to r}g(x),$ and 
use the same convention for any increasing function defined on $(l,r).$ 
The state space of the
diffusion  $Z$ is the interval $(g(l),g(r))$ and
a.s. $\lim_{t\to\zeta(Z)}Z_t=g(r).$ Clearly, letting $y\to r$ in
(\ref{aa10}) it follows that 
\begin{equation}
\label{aa101}
A_{H_r(Y)}=\inf\{t:\, Z_t=g(r)\}\quad {\rm a.s.},
\end{equation}
where both sides in (\ref{aa101}) are either finite or infinite. 
Now we have 

\begin{thm}
\label{!!}
 {\sl For $Y,$ $A,$ $f$ and $g$ as above it holds that $A_{\zeta}$ is
   a.s. finite if and only if for the diffusion $Z$ the boundary point $g(r)$ is an exit
   boundary,  i.e.,
\begin{equation}
\label{aa102}
\int^{g(r)}S^Z(d\al)\int^\al m^Z(d\beta)<\infty,
\end{equation}
where the scale $S^Z$ and the speed $m^Z$ of the
diffusion $Z$ are given by
$$
S^Z(d\al)={\rm e}^{-B^Z(\al)}\,d\al\quad {\rm and}\quad
m^Z(d\beta)=2\,{\rm e}^{B^Z(\beta)}\,d\beta
$$
with 
$$
B^Z(\beta)=2\,\int^\beta G\circ g^{-1}(z)\,dz.
$$
The condition (\ref{aa102}) is equivalent with the condition    
\begin{equation}
\label{aa103}
\int^{r} \left(S^Y(r) - S^Y(v)\right)\, f(v)\, m^Y(dv)<\infty.
\end{equation}
}
\end{thm}
\begin{proof}
As is well known from the standard diffusion theory, a diffusion hits
its exit boundary with positive probability and an exit boundary
cannot be unattainable (see \cite{itomckean74} or
\cite{borodinsalminen02}). This combined with (\ref{aa101}) and the
characterization of an exit boundary (see \cite{borodinsalminen02}
No. II.6 p.14) proves the first claim. It remains to show that
(\ref{aa102}) and (\ref{aa103}) are equivalent. We have
\begin{eqnarray*}
&&
B^Z(\al)=2\,\int^{g^{-1}(\al)} G(u)\, g'(u)\, du\\
&&\hskip1cm
=2\,\int^{g^{-1}(\al)} \left(\frac 12\, \frac{g''(u)}{g'(u)}
+\frac{b(u)}{\sigma^2(u)}\right)\, du\\
&&\hskip1cm
=\log(g'(g^{-1}(\al)))+ B^Y(g^{-1}(\al)).
\end{eqnarray*}
Consequently,
$$
S^Z(d\al)={\rm e}^{-B^Z(\al)}\,d\al= \frac{1}{g'(g^{-1}(\al))}\,
\exp\left(-B^Y(g^{-1}(\al))\right)\,d\al
$$
and
$$ 
m^Z(d\al)=2\,{\rm e}^{B^Z(\al)}\,d\al=2\,g'(g^{-1}(\al))\,
\exp\left(B^Y(g^{-1}(\al))\right)\,d\al.
$$
Substituting first $\al=g(u)$ in the outer integral in (\ref{aa102}) and
after this   $\beta=g(v)$ in the inner integral yield
\begin{eqnarray*}
&&\hskip-.5cm
\int^{g(r)}S^Z(d\al)\int^\al m^Z(d\beta)=2\int^rdu\, {\rm e}^{-B^Y(u)} 
\int^udv\, (g'(v))^2\, {\rm e}^{B^Y(u)}\\
&&\hskip4.2cm
=
2\int^r dv\,(g'(v))^2 \, {\rm e}^{B^Y(v)} \,\int_v^rdu\, {\rm e}^{-B^Y(u)}
\end{eqnarray*}
by Fubini's theorem. Using the expressions given in (\ref{a11}) for the speed and the scale of
$Y$ and the relation (\ref{fg}) between $f$ and $g$ complete the proof. 
\end{proof}

{\bf 4.} It is easy to derive a condition that the mean of $A_\zeta(f)$ is
finite. Indeed, 
\begin{eqnarray}
\label{aa105}
&&
\nonumber
\E_x\left(A_\zeta(f)\right)= 
\int_0^\infty \E_x\left(f(Y_s)\right)\,ds \\
&&\hskip2cm
= \int_l^r G^Y_0(x,y)\, f(y)\,m^Y(dy)<\infty,
\end{eqnarray}
where $G^Y_0$ is the Green kernel of $Y$ w. r. t. $m^Y.$ Under the
assumption (\ref{a10}) we may take for $x\geq y$
$$
G^Y _0(x,y)=S^Y(r) - S^Y(x).
$$
Consequently,  the condition (\ref{aa103}) may be viewed as a part of
the condition (\ref{aa105}). 

{\bf 5.} Since the exit condition (\ref{aa102}) plays a crucial rôle in our approach we discuss here shortly two proofs of this condition, thus making the paper as self-contained as possible. 

Let $Y$ be an arbitrary regular diffusion living 
on the interval $I$ with the end points $l$ and $r.$ The scale function of $Y$ is denoted by $ S$ and the speed measure by $m.$ It is also assumed that the killing measure of $Y$ is identically zero. Recall the definition due to Feller 
\begin{equation}
\label{exit}
r\ {\rm is\ exit}\quad\Leftrightarrow\quad  
\int^{r}S(d\al)\int^\al m(d\beta)<\infty.
\end{equation}
Note that by Fubini's theorem   
$$
\int^{r}S(d\al)\int^\al m(d\beta)=\int^{r}m(d\beta)(S(r)-S(\beta)),
$$
and, hence, $S(r)<\infty$ if $r$ is exit. Moreover, if $r$ is exit then $H_r<\infty$ with positive probability.

\noindent
{\bf 5.1.} We give now some details of the proof of (\ref{exit}) following closely 
Kallenberg \cite{kallenberg97} (see also Breiman \cite{breiman68}). For $l<a<b<r$ let 
$H_{ab}:=\inf\{t : Y_t= a\ {\rm or}\ b\}.$  
Then for  $a<x<b$
\begin{equation}
\label{exit1}
\E_x\left(H_{ab}\right)=\int_a^b \widehat G^Y_0(x,z)\,m(dz),
\end{equation}
where $\widehat G^Y_0$ is the (symmetric) Green kernel of $Y$ killed when it exits $(a,b),$ i.e.,
$$
\widehat G^Y_0(x,z)=\frac{(S(b)-S(x))(S(y)-S(a))}{S(b)-S(a)}\qquad x\geq y.
$$
If $r$ is exit there exists $h>0$ such that $\P_x(H_r<h)>0$ for any fixed $x\in(a,r).$ Using this property it can be deduced (see \cite{kallenberg97} p. 377) that for any $a\in(l,r)$ 
$$
\E_x\left(H_{ar}\right)<\infty,
$$
which, from (\ref{exit1}), is seen to be equivalent with (\ref{exit}). 

\noindent
{\bf 5.2.} Another proof of (\ref{exit}) can be found in Itô and McKean \cite{itomckean74} p. 130. To present also this briefly recall first the formula 
\begin{equation}
\label{for1}
\E_x(\exp(-\lambda\,H_b))=\frac{\psi_\lambda(x)}{\psi_\lambda(b)},
\end{equation}
where $\lambda >0$ and $\psi_\lambda$ is an increasing solution of the generalized differential equation
\begin{equation}
\label{for2}
\frac{d}{dm}\frac{d}{dS} u =\lambda u.
\end{equation}
Letting $b\to r$ in (\ref{for1}) it is seen that 
$$
r\ {\rm is\ exit}\quad\Leftrightarrow\quad  
\lim_{b\to r}\psi_\lambda(b)<\infty.
$$
Let $\psi^+_\lambda$ denote the (right) derivative of $\psi_\lambda$ with respect to $S.$ Since $\psi_\lambda$ is increasing it holds that $\psi^+_\lambda>0.$      
The fact that $\psi_\lambda$ solves (\ref{for2}) yields for $z<r$
$$
\psi^+_\lambda(r)-\psi^+_\lambda(z)=\lambda\int_z^{r} \psi_\lambda(a)\, m(da).
$$
In particular, $\psi^+_\lambda$ is increasing and $\psi^+_\lambda(r)>0.$ Hence, assuming now that $\psi_\lambda(r)<~\infty$ we obtain $S(r)<\infty,$ and, further, 
\begin{eqnarray*}
&&\hskip-.5cm
\lambda\,\psi_\lambda(z) \int^{r}_zS(d\al)\int^\al_z m(d\beta)
\leq \lambda\int^{r}_zS(d\al)\int^\al_z\,\psi_\lambda(\beta) m(d\beta)\\
&&\hskip4.7cm
= \int^{r}_zS(d\al)\left(\psi^+_\lambda(\al)-\psi^+_\lambda(z)\right)\\
&&\hskip4.7cm
= \psi_\lambda(r)-\psi_\lambda(z)  
-\psi^+_\lambda(z)\left(S(r)-S(z)\right)<\infty,
\end{eqnarray*} 
which yields the condition on the right hand side of (\ref{exit}). Assume next that  the condition on the right hand side of (\ref{exit}) holds, and consider for $z<\beta$
$$
0\leq \left(\psi_\lambda(\beta)\right)^{-1}
\left(\psi^+_\lambda(\beta)-\psi^+_\lambda(z)\right)=
\left(\psi_\lambda(\beta)\right)^{-1}\int_z^\beta\psi_\lambda(\al)m(d\al).
$$
Integrating over $\beta$ gives 
\begin{eqnarray*}
&&\hskip-.5cm
\log(\psi_\lambda(r))-
\log(\psi_\lambda(z))
-\psi^+_\lambda(z)\int_z^r 
\left(\psi_\lambda(\beta)\right)^{-1}S(d\beta)\\
&&\hskip3cm
=
\int_z^rS(d\beta)\left(\psi_\lambda(\beta)\right)^{-1}\int_z^\beta\psi_\lambda(\al)m(d\al)\\
&&\hskip3cm
\leq
\int_z^r S(d\beta)\int_z^\beta m(d\al)<\infty,
\end{eqnarray*} 
which implies that $\psi_\lambda(r)<\infty,$ thus completing the proof.

{\bf 6.} As an application of Theorem \ref{!!}, we consider a Bessel process with dimension parameter 
$\delta> 2.$ Let $R$ denote this process. It is well known that $\lim_{t\to\infty} R_t=+\infty$
and that $R$ solves the SDE
$$
dR_t=dW_t+ \frac{\delta-1}{2R_t}\,dt,
$$
where $W$ is a standard Brownian motion. Here the function $B^R$ (cf. (\ref{B1})) takes the form  
$$
B^R(v)=(\delta-1)\log v,
$$
and, consequently,
\begin{eqnarray*}
&&\hskip-.5cm
\int^{\infty}dv\, (g'(v))^2\,{\rm e}^{B^R(v)}\,\int_v^\infty du\, {\rm
  e}^{-B^R(u)}\\
&&\hskip3cm
=\int^{\infty}dv\, (g'(v))^2\,v^{\delta-1}\,\int_v^\infty du\, u^{-\delta+1}\\
&&\hskip3cm
=\int^{\infty}dv\, (g'(v))^2\,v^{\delta-1}\,\frac{1}{\delta-2}\, v^{-\delta+2}
\end{eqnarray*}
leading to 
$$
\int_0^\infty f(R_t)\, dt <\infty \quad \Leftrightarrow\quad
\int^\infty u\, f(u)\, du <+\infty.
$$

Another way to derive this condition is via local times and Jeulin's lemma \cite{jeulin82}. Indeed,
by the occupation time formula and Ray-Knight theorem for the total local times of $R$ (see, e.g.
\cite{yor92c} Theorem 4.1 p. 52)
 we have 
\begin{eqnarray*}
&&
\int_0^\infty f(R_s)\,ds 
\quad{\mathop=^{\rm{(d)}}}\quad
\int_0^\infty f(a)\,\frac{\rho_{a^\gamma}}{\gamma\,a^{\gamma-1}}\, da\\
&&\hskip2.8cm
=\frac 1\gamma\int_0^\infty a\, f(a)\ \frac{\rho_{a^\gamma}}{a^{\gamma}} \,da
\end{eqnarray*}
where $\delta=2+\gamma$ and $\rho$ is a squared 2-dimensional Bessel process. Using 
the scaling property, it is seen that the distribution of the random 
variable $\rho_{a^\gamma}/a^{\gamma}$ does not depend on $a.$ Hence, we obtain by Jeulin's lemma 
that if the function $a\mapsto a\,f(a), \, a>0,$ is locally integrable on $[0,\infty)$ 
then 
\begin{equation}
\label{eq1}
\int_0^\infty f(R_s)\,ds<\infty \quad \Leftrightarrow \quad  \int^\infty a\, f(a)\,da<\infty.
\end{equation}
The same argument allows us to recover the result in \cite{salminenyor03}, that is, 
\begin{equation}
\label{eq2}
\int_0^\infty g(W^{(\mu)}_s)\,ds<\infty \quad \Leftrightarrow \quad\int^\infty  g(x)\,dx<\infty.
\end{equation}
where $g$ is any non-negative locally integrable function and $W^{(\mu)}$ denotes a Brownian motion with drift $\mu>0.$ 
To see this, write $g(x)=f({\rm e}^x)$ 
and use 
Lamperti's representation 
$$
\exp(W^{(\mu)}_s)=R^{(\mu)}_{A^{(\mu)}_s},\quad s\geq 0,
$$
where
$$
A^{(\mu)}_s=\int_0^s du\, \exp(2W^{(\mu)}_u)
$$
and $R^{(\mu)}$ is a Bessel process with dimension $d=2(1+\mu)$ starting from 1, we obtain
(cf. \cite{salminenyor04} Remark 3.3.(3))
$$
\int_0^\infty f(\exp(W^{(\mu)}_s))\,ds=\int_0^\infty 
\left(R^{(\mu)}_u\right)^{-2}\,f(R^{(\mu)}_u)\, du \quad {\rm a.s.},
$$
and, in order to get (\ref{eq2}) it now only remains to use the equivalence (\ref{eq1}). 

We wish to underline the fact that in Theorem \ref{!!} it is assumed that the function $f$ is continuous 
whereas the approach via Jeulin's lemma, which we developed above, demands only local integrability.

\bibliographystyle{plain}
\bibliography{yor1}
\end{document}